\documentclass{amsart}
\usepackage{amsmath,amssymb,amsthm}
\input rlepsf

\def\co{\colon\thinspace}
\newcommand{\Int}{\mbox{\rm Int}}
\newcommand{\N}{\mathbb{N}}
\newcommand{\R}{\mathbb{R}}
\newcommand{\Z}{\mathbb{Z}}
\newcommand{\LL}{\mathbb{L}}
\newcommand{\tb}{{\tt tb}}
\newcommand{\otb}{\overline{\tt tb}}
\newcommand{\lk}{{\tt lk}}
\newcommand{\rot}{{\tt rot}}

\newtheorem{thm}{Theorem}
\newtheorem{lem}[thm]{Lemma}
\newtheorem{prop}[thm]{Proposition}
\newtheorem{cor}[thm]{Corollary}

\theoremstyle{definition}
\newtheorem*{rem}{Remark}
\newtheorem*{remdefn}{Remark/Definition}

\newtheorem*{ex}{Example}
\newtheorem*{ack}{Acknowledgements}

\begin{document}

\title[Legendrian knots and links]{Legendrian knots and links classified by
classical invariants}
\author{Fan Ding}
\address{Department of Mathematics, Peking University,
Beijing 100871, P.~R. China}
\email{dingfan@math.pku.edu.cn}
\author{Hansj\"org Geiges}
\address{Mathematisches Institut, Universit\"at zu K\"oln,
Weyertal 86-90, 50931 K\"oln, Germany}
\email{geiges@math.uni-koeln.de}
\date{}

\begin{abstract}
\noindent It is shown that Legendrian (resp.\ transverse) cable
links in $S^3$ with its standard tight contact structure, i.e.\ links
consisting of an unknot and a cable of that unknot, are classified by
their oriented link type and the classical invariants (Thurston-Bennequin
invariant and rotation number in the Legendrian case, self-linking
number in the transverse case). The analogous result is proved for
torus knots in the $1$--jet space $J^1(S^1)$ with its standard tight
contact structure.

\vspace{.5mm}


\vspace{.5mm}

\end{abstract}

\maketitle

\vspace{2mm}

\section{Introduction}
In \cite{etho01}, Etnyre and Honda explored the methods of convex
surface theory in contact geometry as a tool for the classification of
Legendrian and transverse knots (up to Legendrian or transverse isotopy,
respectively). Their main result was that Legendrian
resp.\ transverse torus knots in $S^3\subset\R^4$ with the standard tight
contact structure
\[ \xi_0=\ker (x_1\, dy_1 -y_1\, dx_1+x_2\, dy_2-y_2\, dx_2)\]
are classified by their topological knot type and the so-called
classical invariants: Thurston-Bennequin invariant $\tb$
and rotation number $\rot$ in
the Legendrian case, self-linking number in the transverse case.
Regarding the former case,
the only previous work in that direction had been the corresponding
result for unknots (in any tight contact $3$--manifold),
due to Eliashberg and Fraser~\cite{elfr98}.
In the transverse case, there are also results due to Eliashberg,
Etnyre, Birman-Wrinkle and Menasco; see~\cite{etho01} for
references.

The aim of the present note is to apply the methods of Etnyre and
Honda to the study of certain Legendrian and transverse links,
where we obtain analogous positive results. The first knot
type whose Legendrian realisations, by contrast, are not determined by the
classical invariants was found by Chekanov~\cite{chek02}; this and further
examples of links with that property are described in
\cite{ng03} and~\cite{ngtr04}.

In order to avoid undue repetition, we assume that the reader is
familiar with the contact geometric concepts discussed in Section~2
of~\cite{etho01}. We shall also have to appeal to a number of results
from convex surface theory as discussed in Section~3 of that paper,
and many of our reasonings are parallel to those employed by Etnyre
and Honda, so we advise the reader to have a copy of their paper at hand
(and better also a copy of Honda's fundamental paper~\cite{hond00}).
Nonetheless, we have tried to give additional details
whenever this seemed to clarify the exposition.
For a brief introduction to convex surface theory one may also
consult~\cite{etny03}; for a survey on knots in contact
geometry see~\cite{etny}.

\section{Cable links in $S^3$}
\label{section:links}
\begin{figure}[h]
\centerline{\relabelbox\small
\epsfxsize 8cm \epsfbox{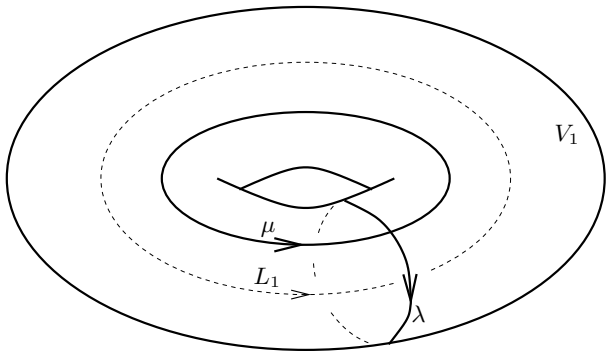}
\extralabel <-2.6cm, 0.4cm> {$\lambda$}
\extralabel <-4.6cm, 1.6cm> {$\mu$}
\extralabel <-4.7cm, 0.9cm> {$L_1$}
\extralabel <-0.7cm, 2.8cm> {$V_1$}
\endrelabelbox}
\caption{Choice of meridian $\mu$ and longitude $\lambda$.}
\label{fig:1}
\end{figure}

Our first result concerns the (oriented)
link $\LL =L_1\sqcup L_2$ in $S^3$, with
$L_1$ a trivial knot and $L_2$ a cable of $L_1$, that is, a torus knot
on the boundary $T$ of a tubular neighbourhood $V_1$ of~$L_1$. Let $\mu$ be
a meridian of the complementary solid torus $V_2:=S^3\setminus\Int (V_1)$,
oriented in the same direction as~$L_1$. Let $\lambda$ a meridian of~$V_1$
(and hence a longitude of~$V_2$),
oriented in such a way that $\mu ,\lambda$ form a positive basis for
$H_1(T)$, with $T$ oriented as the boundary of
$V_2$, see Figure~\ref{fig:1} (where we use the
standard orientation of $\R^3\subset S^3$).

Homologically (on $T$),
$L_2$ is equivalent to $p\mu +q\lambda$ with $p,q$ coprime integers.
Usually we are not going to distinguish between such homologically
equivalent curves and simply say that $L_2$ `is' $p\mu +q\lambda$.
The corresponding link $\LL$ will be referred to as a
{\bf  $(p,q)$--cable link}. (Beware that $L_2$ is a $(p,q)$--torus knot
on~$\partial V_2$, but a $(q,p)$--cable of~$L_1$. Our choice of notation
is more convenient for arguing in analogy with~\cite{etho01}. The only
crucial point, anyway, is the distinction between positive and negative
torus knots, determined by the relative sign of $p$ and~$q$.)
We now consider Legendrian realisations of such links.

\begin{thm}
\label{thm:1}
Two oriented Legendrian cable links in $(S^3,\xi_0)$ are
Legendrian isotopic if and only if their
oriented link types and their classical
invariants agree.
\end{thm}

This theorem will be proved in Section~\ref{section:proof1}.
During the course of that proof we determine --- implicitly, but
completely ---
the range of the classical invariants realised by
cable links. This amounts to a complete classification of such links.

\section{Torus knots in $J^1(S^1)$}
\label{section:knots}
By methods analogous to the proof of Theorem~\ref{thm:1}, we can also
prove a result about Legendrian torus knots in the $1$--jet space
$J^1(S^1)=T^*S^1\times\R$ of the circle with its standard contact
structure~$\xi_1$. This contact structure is defined as
\[ \xi_1=\ker (dz-p\, dq),\]
where $q\in S^1$, the fibre coordinate of $T^*S^1=S^1\times\R$ is denoted
by~$p$, and $z$ is the second $\R$--coordinate. By a {\bf torus knot} in
$J^1(S^1)$ we mean a knot that sits on a torus topologically
isotopic to
\[ \{ (q,p,z)\in J^1(S^1)\co p^2+z^2=1\} .\]
The following proposition must be well-known to some
of the experts, but we have not found it explicitly
stated in the literature.

\begin{prop}
\label{prop:jet}
Let $K_0$ be a Legendrian unknot in $(S^3,\xi_0)$ with $\tb (K_0)=-1$.
Then there is a contactomorphism
\[ f\co(J^1(S^1),\xi_1)\longrightarrow (S^3\setminus K_0,\xi_0) .\]
\end{prop}

\begin{rem}
The condition $\tb (K_0)=-1$, together with the Bennequin inequality
\[ \tb (K)+|\rot (K)|\leq -\chi (\Sigma ),\]
with $\Sigma$ a Seifert surface for the Legendrian knot~$K$ (which in the
case of $K_0$ we may take to be a disc), forces $\rot (K_0)=0$.
Thus, by the theorem of Eliashberg and Fraser~\cite{elfr98}, the
contactomorphism type of $(S^3\setminus K_0,\xi_0)$ does not depend on
the specific choice of~$K_0$.
\end{rem}

Given a Legendrian knot $K$ in $J^1(S^1)$ homotopic to $n$ times
a generator of the fundamental group $\pi_1(J^1(S^1))\cong\Z$, we define
the classical invariants by
\[ \tb (K)=\tb (f(K))+n^2\;\;\mbox{\rm and}\;\;\rot (K)=\rot(f(K)).\]
Proposition~\ref{prop:jet} and the following theorem about Legendrian
torus knots in $J^1(S^1)$
will be proved in Section~\ref{section:jet}. There we also justify
the above formulae for the classical invariants. Our proof is based
on the proof of Theorem~\ref{thm:1}; this again allows one to determine
the range of the classical invariants.

\begin{thm}
\label{thm:2}
Two oriented Legendrian torus knots in $(J^1(S^1),\xi_1)$
are Legendrian isotopic if
and only if their oriented knot types and their classical invariants
agree.
\end{thm}

\section{Proof of Theorem~\ref{thm:1}}
\label{section:proof1}
Let $\LL =L_1\sqcup L_2$ and $\LL'=L_1'\sqcup L_2'$ be two oriented Legendrian
links of the type considered in Theorem~\ref{thm:1} which have the
same link type and classical invariants. In particular, $L_1$ and $L_1'$
are topological unknots with the same classical invariants, hence Legendrian
isotopic by the result of Eliashberg and Fraser~\cite{elfr98}. We may
therefore assume that $L_1'=L_1$.

The Thurston-Bennequin invariant of the unknot can take any negative
integer value; we are going to write $\tb (L_1)=-m$ with $m\in\N$.
This means that we can find an arbitrarily small tubular neighbourhood
$N_1$ of $L_1$ with convex boundary having two dividing curves of slope
$-m$ relative to our chosen $\mu ,\lambda$, i.e.\ dividing curves that
are homologically of the form $\mu -m\lambda$. The rotation number
of the unknot $L_1$ with $\tb (L_1)=-m$ can take any value in the set
\[ \{ -m+1,-m+3,\ldots ,m-3,m-1\}.\]

We shall frequently have to refer to the concept of
convex tori in so-called standard form, so we state the definition
formally, see~\cite[Section~3.2.1]{hond00}.

\begin{remdefn}
The dividing set of a convex torus consists (in a suitable
identification of the torus with $\R^2/\Z^2$) of $2n$ parallel
curves of some slope~$s$. Such a convex torus is said to be in
{\em standard form} if the characteristic foliation consists
of a linear family (of slope $r\neq s$), called the {\em Legendrian ruling},
with singularities along $2n$ {\em Legendrian divides}
parallel to the dividing curves. By a $C^0$--small perturbation
of the torus, the characteristic foliation may be brought into
standard form. In particular, by a perturbation
near the Legendrian divides, the slope $r$ of the Legendrian
ruling can be modified to any value $r'\neq s$.
\end{remdefn}

Write $L_2'$ as $p'\mu+q'\lambda$ with $p',q'$ coprime integers. Since
$\LL$ and $\LL'$ have the same oriented link type, we have
\[ q=\lk (L_1,L_2)=\lk (L_1,L_2') =q'.\]
By changing the orientations of both $L_2$ and $L_2'$, if necessary, we
may assume that $q=q'\geq 0$.

\begin{rem}
By reversing the role of the interior and the exterior of the torus on which
$L_2$ sits, Etnyre and Honda could assume in \cite{etho01} that in
addition $|p|>q$. We are no longer free to make this choice because of the
presence of~$L_1$. This will affect the calculations of certain
maximal Thurston-Bennequin invariants.
\end{rem}

\vspace{1mm}

{\bf Case 1:} $q=q'=0$. In this case, $L_2$ and $L_2'$ are topologically
trivial in $S^3\setminus L_1$, and therefore Legendrian isotopic by
the result of Eliashberg and Fraser.

\vspace{2mm}

{\bf Case 2:} $q=q'=1$. Then $\LL$ and $\LL'$ are Hopf links. The knots
$L_2$ and $L_2'$ have tubular neighbourhoods $N_2$ and $N_2'$, respectively,
with convex boundary having two dividing curves of the form
$\lambda +\tb (L_2) \mu$, i.e.\ of slope $1/\tb (L_2)$. This allows one
to find a contactomorphism $\phi_2\co N_2\rightarrow N_2'$ sending
$L_2$ to~$L_2'$ and preserving the choice of longitude~$\lambda$
(and of course the meridian~$\mu$).

We may choose a standard tubular neighbourhood $N_1$ of $L_1$
such that $N_1\cap N_2=N_1\cap N_2'=\emptyset$. Set
\[ M= S^3\setminus \Int (N_1\cup N_2)\cong T^2\times [0,1]=
\mu\times\lambda\times [0,1].\]
The Legendrian rulings on the
two convex boundary components of $M$ give a section $\sigma$ of
$\xi_0$ over $\partial M$. In the local model
\[ \cos\theta\, dx-\sin\theta\, dy=0\]
for a neighbourhood of a Legendrian curve $\{ x=y=0\}$, this simply
corresponds to the vector field~$\partial_{\theta}$. In particular,
$\sigma$ extends over $N_1$ and $N_2$ so as to coincide with the
section of $\xi_0$ defined by the oriented Legendrian knots $L_1$ and~$L_2$.
It follows that the relative Euler class $e(\xi_0 ,\sigma )
\in H^2(M,\partial M)$
evaluates (up to fixing signs)
to $\rot (L_1)$ and $\rot (L_2)$ on the annuli $\mu\times
[0,1]$ and $\lambda\times [0,1]$, respectively (cf.~\cite[p.~354]{hond00}).
This determines $e(\xi_0,\sigma )$ completely.

The contact structure $\xi_0$ restricted to the thickened torus $M$ is
minimally twisting in the sense of~\cite[p.~313]{hond00} --- this corresponds
to the $2\pi$--torsion in the sense of~\cite[Defn.~1.2]{giro00} being zero
---, for otherwise we would find an overtwisted disc in~$(S^3,\xi_0)$.
Therefore, from the classification of such contact structures on a thickened
torus~\cite[Prop.~4.22 resp.\ Prop.~4.9]{hond00},
we deduce that the identity map on $N_1$
and the contactomorphism $\phi_2\co N_2\rightarrow N_2'$ extend to
a contactomorphism of $(S^3,\xi_0)$. By a result of
Eliashberg~\cite[Cor.~2.4.3]{elia92}, any such contactomorphism is contact
isotopic to the identity. It follows that $\LL$ and $\LL'$ are
Legendrian isotopic.

\begin{rem}
In the above argument, where we appeal to Honda's classification
of minimally twisting contact structures on a thickened torus,
his Proposition 4.22 is used when the slopes of the dividing
curves on the two boundary tori are different. Here the
different contact structures are distinguished by
the relative Euler class, which we have shown to be
determined by $\rot (L_1)$ and $\rot (L_2)$. Honda's Proposition 4.9 is
used in the nonrotative case, when the two boundary slopes are equal.
This is the case if and only if $\tb (L_1)=\tb (L_2)=-1$
(and $\rot (L_1)=\rot (L_2)=0$).
\end{rem}

\begin{rem}
It is easy to see in the front projection picture for Legendrian knots
that in Cases 1 and 2 each of the two Legendrian unknots $L_1,L_2$
can realise any combination of $\tb$ and $\rot$ allowed for a single unknot.
\end{rem}

\vspace{1mm}

{\bf Case 3:} $q=q'\geq 2$. In this case we must have $p,p'\neq 0$, and we
claim that $p=p'$. Indeed, if $p\neq\pm 1$, this follows from the
classification of torus knots~\cite[Thm.~3.29]{buzi85}. The case $p=\pm 1$
is equivalent to $L_2$ (and hence $L_2'$) being a trivial knot, so we
deduce $p'=\pm 1$, but potentially $p,p'$ of opposite sign. We see that
$p=p'$ in this case as well by considering the $q$--fold covering
of $S^3$ by itself, branched along~$L_1$. The inverse image of $L_2$ under
this covering has $q$ components, each representing $p\mu +\lambda$.
Any two of these components have linking number~$p$. For $L_2'$ the
situation is analogous, and since the link types of $\LL$ and $\LL'$ agree
we conclude that $p=p'$.

We now have to distinguish two subcases, distinguished by the sign of~$p$.
It will be understood from now on that $q=q'$ and $p=p'$.

\vspace{2mm}

{\bf Case 3 (a):} $q\geq 2$ and $p>0$. The maximal Thurston-Bennequin
invariant of a Legendrian $(p,q)$--torus knot $L_2$ with $p,q>0$ (including the
case $p=1$, i.e.\ the trivial knot) is $pq-p-q$, see~\cite[Thm.~4.1]{etho01}.
The rotation number of a Legendrian torus knot realising this maximal
Thurston-Bennequin invariant equals~$0$; this follows immediately
from the Thurston-Bennequin inequality and the fact that the (positive)
$(p,q)$--torus knot has a Seifert surface of Euler characteristic
$p+q-pq$, cf.~\cite[Cor.~4.11]{buzi85} and~\cite{miln68}.
An explicit example of such a Legendrian torus knot
is shown in Figure~8 (top) of~\cite{etho01}; this front projection picture
is valid without the restriction $p>q$. We can add a meridional Legendrian
unknot~$L_1$ (with any possible combination
of $\tb$ and $\rot$) to that picture so as to obtain a Legendrian
realisation of our cable link $\LL$, so the maximal Thurston-Bennequin
invariant of $L_2\subset\LL$ has that same value $pq-p-q$.

\begin{lem}
\label{lem:isotopic}
If $\tb (L_2)=\tb (L_2')=pq-p-q$, then $\LL$ and $\LL'$ are
Legendrian isotopic.
\end{lem}

\begin{proof}
The argument is largely parallel to the proof of the corresponding Lem\-ma~4.7
of~\cite{etho01}. Using the reasoning of Etnyre and Honda, we obtain
two Heegaard splittings $V_1\cup V_2$ and $V_1'\cup V_2'$ of $S^3$ into
two solid tori with the following properties:
\begin{itemize}
\item $T:=V_1\cap V_2$ is a convex torus on which $L_2$ sits,
\item $L_1$ is contained in the interior of $V_1$,
\item $V_2$ is a standard neighbourhood of a Legendrian unknot with
$\tb =-1$ (and $\rot =0$);
\end{itemize}
analogously for the primed spaces.
We may assume in addition that $N_1\subset V_1\cap V_1'$.

We briefly elaborate on that third point (the characterisation of~$V_2$),
since variants of this argument will be used below. The boundary torus
$T=\partial V_2$ may be assumed in standard form, since the existence
of additional intersection points of $L_2$ with the dividing set would
allow a stabilisation of $L_2$, which is impossible for $\tb (L_2)$
being maximal. Moreover, there is a unique tight contact structure
on a solid torus with a fixed convex boundary of slope~$-1$, so one finds a
contactomorphism $V_2\rightarrow V_2'$.

As in the argument for Case~2 (with $V_2$ taking the role
of~$N_2$), $\rot (L_1)$ determines the contact
structure on the thickened torus $S^3\setminus (N_1\cup V_2)$. This then
yields a contact isotopy that moves $L_1$ to itself and $T$ to $T'$. The
ruling curves on $T=T'$ then define a Legendrian isotopy of $L_2$ to
$L_2'$ that extends to a contact isotopy fixing~$L_1$. (For a proof of
the Legendrian isotopy extension theorem see~\cite{geig}.)
\end{proof}

The proof of Case~3~(a) is now completed by showing as in
Lemma~4.8 of~\cite{etho01} that if $\tb (L_2)$ is not maximal, then it
is a stabilisation of a Legendrian knot $L_2^*\subset S^3\setminus L_1$
(with maximal Thurston-Bennequin invariant). By the last lemma, this determines
the link $L_1\sqcup L_2^*$ up to Legendrian isotopy.
The number of positive and the number of negative stabilisations
are determined by $\tb (L_2)$ and $\rot (L_2)$.
There is only one small point where one
has to take extra care:

Again we are dealing with a Heegaard splitting
$S^3=V_1\cup V_2$ into solid tori, with $L_1\subset\Int (V_1)$
and $T=V_1\cap V_2$
a torus on which $L_2$ sits. If the slope of $T$, when regarded as the boundary
of~$V_2$, is less than or equal to $-1$, the argument of Etnyre and Honda
applies without any changes. If the slope of $T$ is less than $-1$ when
regarded as the boundary of~$V_1$, we need to ensure that we can find
a convex torus $T'$ with two dividing curves of slope $-1$ not simply
in $V_1$, but in $V_1\setminus N_1$. Since $N_1$ may be taken to be a standard
neighbourhood of $L_1$ with two dividing curves on $\partial N_1$ of
slope $1/\tb (L_1)\geq -1$ (relative to meridian and longitude adapted to
$N_1$!), and slopes decrease as we move outwards, this is indeed possible:
Lemma~3.15 of \cite{etho01} guarantees that we can find a torus $T_0$
of slope~$-1$, but possibly with more than two dividing curves.
Theorem~2.2 (1) of \cite{hond00} then allows us to split off a nonrotative
collar $T_0\times [0,1]$ of $T_0\equiv T_0\times\{ 0\}$
on the toric annulus between $T_0$ and $\partial V_1$, such that the torus
$T_1:=T_0\times\{1\}$ in the interior of this toric annulus has the same
slope as $T_0$, but only two dividing curves.
(The existence statement of that theorem is not affected by the erratum.)
This is all the information that is necessary to ensure that $L_2$
can be destabilised.

\vspace{2mm}

{\bf Case 3 (b):} $q\geq 2$ and $p<0$. Again the first task we have to
deal with is to determine the maximal Thurston-Bennequin invariant of~$L_2$.

If $p<-1$, then $L_2$ is a nontrivial negative torus knot, in which case
Lemma~4.9 of~\cite{etho01} tells us that $\tb (L_2)\leq pq$. We want to obtain
the same estimate in the case $p=-1$.

\begin{lem}
\label{lem:tb-estimate}
If $p=-1$, then
$\tb (L_1)\leq -q$ or $\tb (L_2)\leq -q$.
\end{lem}

\begin{proof}
Suppose that $\tb (L_1)>-q$. Any $(\tb (L_1)+q -1)$--fold
stabilisation $\tilde{L}_1$ of
$L_1$ has $\tb (\tilde{L}_1)=-q+1$.
Perform Legendrian surgery on such a stabilised knot~$\tilde{L}_1$,
i.e.\ contact
$(-1)$--surgery in the sense of~\cite{dige04}. The result of this surgery
will be a fillable and hence tight contact structure on some lens space.
The framing of this surgery relative to the surface framing of $\tilde{L}_1$
is $\tb (\tilde{L}_1)-1 =-q$. With our definition of $\mu$ and $\lambda$,
this means that the curve $\mu -q\lambda$ (and hence~$L_2$) bounds
a disc in the surgered manifold.

If we had $\tb (L_2)\geq -q$, we could stabilise $L_2$ so as to obtain
a Legendrian knot $\tilde{L}_2$ topologically isotopic to $L_2$
and with $\tb (\tilde{L}_2)=-q$. But $-q$ is also the linking of
$L_2$ with a parallel copy of it on the torus on which it sits.
This implies that
in the surgered manifold the contact framing of $\tilde{L}_2$ and its
surface framing coincide, making it the boundary of an overtwisted
disc. This contradicts the tightness of the surgered manifold.

Thus, if $\tb (L_1)>-q$, we must in fact have $\tb (L_2)\leq -q-1$
(but we only need the weaker statement of the lemma).
\end{proof}

Therefore --- still in the case $p=-1$ ---, if $\tb (L_2)\leq -q$, we have the
desired estimate $\tb (L_2)\leq pq$ in this case as well. If, however,
$\tb (L_2)>-q$, then the lemma tells us that $\tb (L_1)\leq -q$.
Since $p=-1$, we may regard $L_1$ as a cable of the unknot~$L_2$.
In terms of the corresponding choice of meridian $\mu_2$ and longitude
$\lambda_2$ (of the complement of a tubular neighbourhood $V_2$ of~$L_2$!)
the class of $L_1$ on $\partial V_2$ is $-\mu_2+q\lambda_2$.

Thus, by interchanging the roles of $L_1$ and $L_2$, if necessary,
we may always assume that $\tb (L_2)\leq pq$. In the next lemma
we determine the actual maximal Thurston-Bennequin invariant $\otb (L_2)$
of $L_2$ in a Legendrian realisation of our cable link. Recall that
we wrote $m=-\tb (L_1)$.

\begin{lem}
\label{lem:maxtb}
We have $\otb (L_2)=pq-\max (mp+q,0)$.
\end{lem}

\begin{rem}
If $|p|>q$, this yields $\otb (L_2)=pq$ as in~\cite{etho01}.
\end{rem}

\begin{proof}[Proof of Lemma~\ref{lem:maxtb}]
First we are going to find a Legendrian cable link $\LL$ where $L_2$
realises the claimed maximal Thurston-Bennequin invariant. We may assume
that $\partial V_1$ is convex with two dividing curves of the
form $\mu -m\lambda$, i.e.\ curves of slope $-m$ when seen from
the exterior solid torus $S^3\setminus\Int (V_1)$, and characteristic foliation
in standard form.

(i) In the case that $mp+q\leq 0$, that is, $q/p\geq -m$,
there is (by Lemma~3.16 of~\cite{etho01}) a convex torus in
$S^3\setminus V_1$, parallel to $\partial V_1$, with two dividing curves of
slope $q/p$. Take $L_2$ to be one of its two Legendrian divides,
which have the same slope~$q/p$,
so that $L_2$ is of the form $p\mu +q\lambda$. For a Legendrian
divide, the contact framing coincides with the framing it inherits from
the torus. This implies $\tb (L_2)=pq$.

(i') If $mp+q>0$, then in particular $q/p\neq -m$, so we may assume
that the Legendrian ruling curves on $\partial V_1$ have slope $q/p$. 
Let $L_2$ be such a ruling curve. We claim that $\tb (L_2)=pq-(mp+q)$.
For that we have to appeal to the (second part of the)
following result of Kanda~\cite{kand98},
cf.~\cite[Thm.~3.4]{etho01}:

\begin{thm}[Kanda]
\label{thm:Kanda}
If $\gamma$ is a Legendrian curve in a surface $\Sigma$, then $\Sigma$
may be isotoped relative to $\gamma$ so that it is convex if and only if
the twisting $t_{\Sigma}(\gamma )$ of the contact planes along $\gamma$
relative to the framing induced by $\Sigma$ satisfies $t_{\Sigma}(\gamma )
\leq 0$. If $\Sigma$ is convex, then
\[ t_{\Sigma}(\gamma )=-\frac{1}{2}\# (\gamma\cap\Gamma ),\]
where $\# (\gamma\cap\Gamma )$ denotes the number of intersection points
of $\gamma$ with the dividing set $\Gamma$ of~$\Sigma$.
\hfill\qed
\end{thm}

This implies that if $\gamma'$ is a parallel curve to $\gamma$ in~$\Sigma$,
then
\begin{eqnarray*}
\tb (\gamma ) & = & \lk (\gamma ,\gamma')+t_{\Sigma}(\gamma )\\
              & = & \lk (\gamma ,\gamma')-\frac{1}{2}\# (\gamma\cap\Gamma ).
\end{eqnarray*}
In our case this yields
\[ \tb (L_2)=pq-\left|\begin{array}{rr}p&-1\\q&m\end{array}\right|
= pq-(mp+q).\]

(ii) It remains to show that the values for $\tb (L_2)$ that we have
found in (i) and (i') are the maximal possible. Thus, let $\LL =L_1\sqcup
L_2$ be a Legendrian realisation of our $(p,q)$--cable link.
Let $T$ be a standardly embedded torus in $S^3\setminus L_1$ on
which $L_2$ sits. The linking of $L_2$ with a push-off along $T$
equals~$pq$. As we observed before this lemma, we have
$\tb (L_2)\leq pq$.
It follows that $t_T(L_2)\leq 0$, so by Kanda's theorem we can make
$T$ convex without moving~$L_2$. Furthermore, we only need to deal with
the case $mp+q>0$.

We may choose a small standard neighbourhood $N_1$ of $L_1$ such that
$T\subset S^3\setminus N_1$.
The slope of the dividing curves of $T$ must be negative and greater
than or equal to~$-m$ (the slope of~$\partial N_1$),
otherwise we would find (again by Lemma~3.15
of~\cite{etho01}) a torus between $T$ and $\partial N_1$ of slope~$0$,
giving rise to an overtwisted disc (and hence a contradiction). Write the
slope of $T$ as $-r/s\geq -m$ with $r,s\geq 1$. Observe that, since $p<0$,
\[ rp+sq\geq s(mp+q)>0,\]
so the algebraic intersection number of $L_2$ with a dividing curve
equals $\left|\begin{array}{rr}p&-s\\q&r\end{array}\right|$.
The geometric intersection number has to be at least as big as that, so
Kanda's formula yields, assuming that there are $2n$ dividing curves,
\[ \tb (L_2)\leq pq-n\left|\begin{array}{rr}p&-s\\q&r\end{array}\right|
\leq pq-(mp+q).\]
This concludes the proof of Lemma~\ref{lem:maxtb}.
\end{proof}

Classifying Legendrian links realising this maximal Thurston-Bennequin
invariant is much more involved in Case 3~(b) than it was in Case 3~(a),
since here the rotation number can take on different values
(at least if $mp+q\leq 0$). This also implies
that stabilisations of different such links can yield Legendrian isotopic
links.
Before we address these issues, we show as in Case 3~(a) that links with
non-maximal Thurston-Bennequin invariant can be destabilised. 

\begin{lem}
\label{lem:destabilise}
Let $L_1\sqcup L_2$ be a Legendrian realisation of a $(p,q)$--cable
link. If $\tb (L_2)<pq-\max (mp+q,0)$, then there is a link $L_1\sqcup
L_2^*$ such that $\tb (L_2^*)>\tb (L_2)$ and $L_2$ is a stabilisation
of $L_2^*$ in $S^3\setminus L_1$.
\end{lem}

\begin{proof}
Since $\tb (L_2)<pq$, Kanda's theorem allows us to assume that
$L_2$ lies on a convex, standardly embedded torus $T$ in $S^3\setminus L_1$.
Let $-r/s$ be the slope of (the dividing curves of)~$T$.
If $-r/s=q/p$, then the algebraic intersection number of $L_2$ with the
dividing set $\Gamma$ of $T$ is zero, but the actual
geometric intersection number $\# (L_2\cap\Gamma )$ is positive by
Kanda's formula. In this case, $L_2$ can be destabilised
(see~\cite[p.~85]{etho01}).

We may therefore assume that $-r/s\neq q/p$.
Let $V_1\cup V_2$ be the Heegaard
splitting associated with $T$, where $N_1\subset \Int (V_1)$ may be taken
as a standard neighbourhood with convex boundary $\partial N_1$ having
two dividing curves of slope~$-m$ and in standard form.

(i) If $mp+q\leq 0$, then $q/p\geq -m$. By the same reasoning as above
we can find a convex torus $T_{q/p}$ of slope $q/p$ either in~$V_2$
(if $q/p>-r/s$),
or in $V_1\setminus N_1$ (if $-r/s>q/p\geq -m$). Connect $L_2$ with a 
Legendrian divide $\gamma$ on $T_{q/p}$ by an annulus~$A$. Since $t_A
(\gamma )= t_{T_{q/p}} (\gamma )=0$ (the surface framing
coincides with the contact framing), and $t_A(L_2)=t_T(L_2)=\tb (L_2)-pq<0$,
Kanda's theorem allows us to assume that $A$ is convex. (In fact, it
is enough to know that $\otb\leq pq$ for a $(p,q)$--torus knot
in order to deal with both $\gamma$ and $L_2$ at one stroke.)
Now the dividing curves
on $A$ will intersect $L_2$ more often than~$\gamma$. This gives
rise to a bypass and hence a destabilisation of~$L_2$. 

(ii) If $mp+q>0$, then $-r/s\neq -m$ (hence $-r/s>-m$, cf.\ part (ii)
of the proof of Lemma~\ref{lem:maxtb}) or the number $2n$ of dividing curves
on $T$ is greater than~$2$. Otherwise, by Kanda's formula, we would have
\[ \tb (L_2)=pq-\left|\begin{array}{rr}p&-1\\q&m\end{array}
\right|=pq-(mp+q).\]

Since $q/p\neq -m$, we may assume that the Legendrian ruling curves on
$\partial N_1$ have slope~$q/p$. Let $A$ be an annulus connecting $L_2$
with such a ruling curve. As in (i), we may take $A$ to be convex. The number
of intersection points of $L_2$ with the dividing curves of $A$
is equal to
\[ 2n\left|\begin{array}{rr}p&-s\\q&r\end{array}\right|=2n (rp+sq)
\geq 2ns(mp+q)\geq 2(mp+q)=2\left|\begin{array}{rr}p&-1\\q&m\end{array}
\right|,\]
that last number being equal to the number of intersection points of
a Legendrian ruling curve on $\partial N_1$ with the dividing curves of~$A$.
At least one of the two inequalities is strict.
As in (i), this allows one to perform a destabilisation of~$L_2$.
\end{proof}

It turns out that the case with $mp+q\leq 0$ is analogous to
the case of negative torus knots in~\cite{etho01}, whereas the case
with $mp+q>0$ can be treated much more simply by methods analogous to the
case of positive torus knots. We begin with that latter case.

\vspace{2mm}

{\bf Case 3 (b1):} $q\geq 2$, $p<0$, and $mp+q>0$. Here, analogous to
Lemma~\ref{lem:isotopic}, the following lemma says that the links
with $L_2$ realising the maximal Thurston-Bennequin invariant
are Legendrian isotopic. Together with Lemma~\ref{lem:destabilise}
this finishes the proof of Theorem~\ref{thm:1} in the present case.

\begin{lem}
If $\tb (L_2)=\tb (L_2')=pq-(mp+q)$, then $\LL$ and $\LL'$ are
Legendrian isotopic.
\end{lem}

\begin{proof}
Let $T$ be a standardly embedded torus in $S^3\setminus L_1$ on which
$L_2$ sits. By Kanda's theorem we can make $T$ convex without moving~$L_2$.
Moreover, the fact that $\tb (L_2)$ is maximal allows us to
assume that $T$ is in standard form.
In the computations in part (ii) of the proof of
Lemma~\ref{lem:maxtb} we have equality everywhere (i.e.\ maximal
Thurston-Bennequin invariant) if and only if the slope $-r/s$ of $T$
equals $-m$ and there are $2n=2$ dividing curves.

The argument is now analogous to the proof of Lemma~\ref{lem:isotopic},
but with one small added complication regarding the classification
of tight contact structures on~$V_2$ (with notation as in
Lemma~\ref{lem:isotopic}). This classification is essentially
determined by the dividing set on $\partial V_2$ and not
affected by assuming $\partial V_2$ to be standard with Legendrian ruling
of slope~$0$ (i.e.\ with Legendrian meridians), see~\cite[Prop.~4.2]{hond00}.
In that case, contact structures on $V_2$ are classified by
the rotation number of the meridian, see~\cite[Prop.~4.23]{hond00}.

Since $T$ (and likewise $T'$)
and $\partial N_1$ have the same slope and two dividing curves each,
we are in the nonrotative case dealt with in~\cite[Prop.~4.9]{hond00}.
According to that proposition, there is a unique contact structure
(up to diffeomorphism) on the thickened torus $S^3\setminus (N_1\cup V_2)$.
Moreover, this implies that the rotation number of a meridian
of $V_2$ equals $\rot (L_1)$. By what we said above, this allows us
to find a contactomorphism $V_2\rightarrow V_2'$.
The proof now concludes like that of Lemma~\ref{lem:isotopic}.
\end{proof}

From this lemma we immediately conclude that the rotation number of a
Legendrian knot $L_2$ realising this maximal Thurston-Bennequin invariant
$pq-(mp+q)$ (in a cable link) must be fixed by the data for~$L_1$.
The following lemma establishes that rotation number.

\begin{lem}
\label{lem:maxtb-rot}
If $\tb (L_2)=pq-(mp+q)$, then $\rot (L_2)=p\cdot\rot (L_1)$.
\end{lem}

\begin{rem}
The method of proof used here is related to that introduced by Etnyre
and Honda for negative torus knots~\cite[pp.~88--93]{etho01}, and one that we
shall also have to use in greater generality below for the remaining
case $mp+q\leq 0$.
\end{rem}

\begin{proof}[Proof of Lemma~\ref{lem:maxtb-rot}]
As seen in part (i') of the proof of Lemma~\ref{lem:maxtb}, $L_2$ can
be realised as a Legendrian ruling curve of slope
$q/p$ on a standard convex torus of slope~$-m$. The standard
model for a convex torus of slope $-m$ (of the dividing curves)
and Legendrian ruling of slope $\infty$ is given by $\partial N_1$, where
\[ N_1=\{ (x,y,\theta )\in\R^2\times S^1\co x^2+y^2\leq 1\}\]
with contact structure
\[ \cos (m\theta )\, dx-\sin (m\theta )\, dy =0\]
on $\R^2\times S^1$. (Remember that we are measuring slopes
with respect to the complementary solid torus $S^3\setminus N_1$.)
The dividing set of $\partial N_1$ is
\[ \Gamma =\{ (\pm\sin (m\theta ),\pm\cos (m\theta ),\theta )\co \theta
\in S^1\};\]
the Legendrian ruling curves are the $\theta$--coordinate lines on~$\partial
N_1$.

Now regard $N_1$ as a neighbourhood of $L_1\subset (S^3,\xi_0)$,
with $L_1\subset N_1$ given by $\{ x=y=0\}$ (and oriented by
$\partial_{\theta}$).
Strictly speaking, we may have to take $N_1$ of smaller radius in the
$xy$--plane, but this does not affect the following homotopical arguments.
The rotation number $\rot (L_1)$ measures the rotation of the
tangent vector field $\partial_{\theta}$ along
$L_1$ relative to a trivialisation of~$\xi_0$. The same counting of
rotations of $\partial_{\theta}$ relative to a trivialisation of $\xi_0$
is possible along any closed curve, even if $\partial_{\theta}$ is
not tangent to that curve. This associates a 
`generalised' rotation
number with any (not necessarily Legendrian) closed curve
in~$N_1$. In particular, the curve $p\mu +q\lambda$ of
slope $q/p$ will have generalised rotation number~$p\cdot\rot (L_1)$, since
it is homotopic in $N_1$ to the curve $pL_1$.

It is now possible to make a $C^0$--small perturbation of $\partial N_1$
(through tori in standard form) such
that it becomes a torus in standard form with Legendrian rulings of
slope $q/p$ (and divides still of slope~$-m$), cf.~\cite[Cor.~3.6]{hond00}.
Let $v$ be a non-zero section of $\xi$ along the curve $p\mu +q\lambda$
defined by the Legendrian ruling curves (on each perturbed torus).
Thus, we start out with $v=\partial_{\theta}$ and end with a tangent
vector field to a Legendrian $(p,q)$--ruling curve~$L_2$. The rotation number
$\rot (L_2)$ is the rotation of $v$ relative to a trivialisation of~$\xi_0$.
By the described deformation argument, this equals the rotation of
$v=\partial_{\theta}$ along the curve $p\mu +q\lambda$ on~$\partial N_1$,
that is, $p\cdot\rot (L_1)$.
\end{proof}

\begin{rem}
The preceding proof shows in particular how to realise a Legendrian link
with any allowable classical invariants for the unknot $L_1$ and
$\tb (L_2)$ equal to the maximum~$\otb (L_2)=pq-(mp+q)$;
according to the lemma, $\rot (L_2)$ is then determined.
\end{rem}

\begin{ex}
Here is an explicit realisation of the cable link with $\tb (L_2)=
\otb (L_2)$ in the case that $p=-1$ (and the other conditions
of Case~3~(b1), so in particular $q>m$). Regard $\LL$ as a link
in $\R^3$ with its standard tight contact structure
\[ \xi_0=\ker (dz+x\, dy);\]
this space is contactomorphic to the complement of a point in
$(S^3,\xi_0)$, cf.~\cite{geig}. Re\-present the oriented knot
$L_1$ by its front projection to the
$yz$--plane. Define the Legendrian knot
$L_2$ by the following front projection picture:
Begin with a Legendrian push-off $L_2^0$ of $L_1$, given
by a parallel copy of the front projection, but with
{\em reversed} orientation (to ensure $p=-1$). Then $\tb (L_2^0)=
\tb (L_1)$ and $\rot (L_2^0)=-\rot (L_1)$. Moreover, we have
$\lk (L_1,L_2^0)=-\tb (L_1)=m$. Next, add $q-m$ meridional loops
around $L_1$ to $L_2^0$ as shown in Figure~\ref{figure:loops};
take $L_2$ to be the resulting knot.

\begin{figure}[h]
\centerline{\relabelbox\small
\epsfxsize 8cm \epsfbox{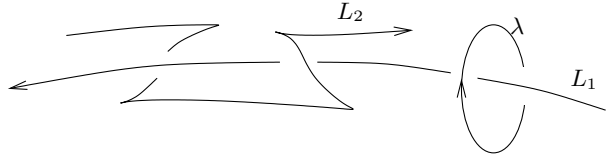}
\extralabel <-0.5cm, 0.9cm> {$L_1$}
\extralabel <-1.3cm, 1.6cm> {$\lambda$}
\extralabel <-3.6cm, 1.8cm> {$L_2$}
\endrelabelbox}
\caption{Adding a meridional loop.}
\label{figure:loops}
\end{figure}

Then $\lk (L_1,L_2)=m+(q-m)=q$, so $L_2$ represents $-\mu +q\lambda$.
For each of the $q-m$ meridional loops we add four cusps to the front
projection (but leave its writhe unchanged), so
\[ \tb (L_2)=\tb (L_2^0) -2(q-m)=-m-2(q-m)=pq-(mp+q)
=\otb (L_2).\]
Since we have added as many up-cusps as down-cusps by introducing extra
meridional loops, we have $\rot (L_2)=\rot (L_2^0)=-\rot (L_1)$.
Here we have used the well-known formulae for computing $\tb$ and $\rot$
from the front projection, cf.~\cite{etny03}. In the next section
these formulae are stated for $J^1(S^1)$, but they hold
equally in~$(\R^3,\xi_0)$ .
\end{ex}

\vspace{1mm}

{\bf Case 3 (b2):} $q\geq 2$, $p<0$, and $mp+q\leq 0$.
Assume that $\LL$ is a Legendrian realisation of our $(p,q)$--cable link
with $\tb (L_2)=pq$ taking on the maximal possible value.
By the now familiar argument involving Kanda's
theorem we may assume that $L_2$ sits
on a convex torus $T$ standardly embedded in $S^3\setminus L_1$.
Since $\tb (L_2)$ equals the linking of $L_2$ with a push-off along~$T$,
Kanda's formula implies that the slope of $T$ must be~$q/p$.
Moreover, it may be assumed that $T$ is in
standard form, the number of Legendrian divides
is two, and that $L_2$ is one of them (cf.\ the end of the proof
of Lemma~4.8 in~\cite{etho01} and appeal to Giroux flexibility
rel boundary~\cite{giro91}, cf.~\cite[Thm.~3.4]{hond00}).

Our first task is to show that the other classical invariants,
i.e.\ $\tb (L_1)$, $\rot (L_1)$ and $\rot (L_2)$, then determine $\LL$
up to Legendrian isotopy. Moreover, one has to determine the range of
$\rot (L_2)$ and, finally, to discuss the equivalence
of Legendrian links $\LL$, $\LL'$ with equal invariants that
destabilise to links with different invariants.

In order to deal with these issues, we need to recall an invariant of
homology classes of curves on a convex torus $T$ in standard form, as
defined in~\cite[p.~88]{etho01}. Let $v_0$ be a globally non-zero section of
$\xi_0$, and $v_T$ a section of $\xi |_T$ that is tangent to the
Legendrian ruling curves and
transverse to the Legendrian divides.
(In \cite{etho01}, $v_T$ is chosen transverse to the ruling
and tangent to the Legendrian divides, but up to homotopy this is the same
as our choice.) For $\gamma$ a closed, oriented curve
on $T$, define $f_T(\gamma )$ as the rotation of $v_T$ relative to $v_0$
along~$\gamma$.
This function $f_T$ has the following properties, which are not difficult
to check (for the second one cf.\ the argument in the proof of
Lemma~\ref{lem:maxtb-rot} above).

\begin{itemize}
\item The function $f_T$ is well-defined and linear on homology classes.
\item The function $f_T$ is invariant under isotopies of $T$ through
convex tori in standard form.
\item If $\gamma$ is a Legendrian ruling curve or a Legendrian divide,
then $f_T(\gamma )=\rot (\gamma )$.
\end{itemize}

In particular, we have
\[ \rot (L_2)=pf_T(\mu )+qf_T(\lambda).\]

The goal now is to show that $\rot (L_2)$ actually determines
$f_T(\mu )$ and $f_T(\lambda )$, and that this is the key to
completing the classification.
Since Etnyre and Honda may assume without loss of generality that
$|p|>q$, but we --- unfortunately --- may not, we have to make
one last division into three subcases.

\vspace{2mm}

{\bf Case 3 (b2-i):} $q\geq 2$, $p<0$, $mp+q\leq 0$, and $|p|>q$.
In this case, the arguments of~\cite{etho01} apply, and we only need
to add a few details.
Etnyre and Honda derive various properties of the function~$f_T$,
notably the range of the values $f_T(\mu )$ and $f_T(\lambda )$.
From there one concludes that the rotation number $\rot (L_2)$ does indeed
determine both $f_T(\mu )$ and $f_T(\lambda )$. Moreover, the
possible values of $\rot (L_2)$ (for $\tb (L_2)=pq$ being maximal)
are shown to lie in
\[ \left\{ \pm (p+(2l+1)q)\co l\in\Z ,\; 0\leq l< -\frac{p+q}{q}\right\} .\]
All these rotation numbers are realised by a Legendrian
knot $L_2$ forming part of a cable link $\LL$: in Figure~8 (bottom)
of~\cite{etho01} one
finds an explicit realisation of $L_2$, to which a meridional loop
$L_1$ can be added with $\tb (L_1)=-m$ and $\rot (L_1)$
equal to one of the allowable values $-m+1, -m+3,\ldots ,m-3,m-1$.

\begin{lem}
\label{lem:3b2i-isotopy}
Let $\LL$ and $\LL'$ be $(p,q)$--cable links with
$\tb (L_2)=\tb (L_2')=pq$. Then
$\LL$ and $\LL'$ are Legendrian isotopic if and only if $\rot (L_2)
=\rot (L_2')$.
\end{lem}

\begin{proof}
We split $S^3$ into a solid torus $V_2$ (with boundary~$T$, and
$L_1$ in the complementary solid torus), a standard
neighbourhood $N_1$ of~$L_1$, and a thickened torus $M$ with boundaries
$T$ and $\partial N_1$; similarly we define $V_2'$ and~$M'$.

The rotation number $\rot (L_2)$ determines $f_T(\mu)$ (which
equals the rotation number of $\mu$ when the Legendrian ruling is
made horizontal), and hence the contact structure on
$V_2$ resp.\ $V_2'$ by the classification of tight
contact structures on solid tori~\cite[Prop.~4.23]{hond00}. Thus,
we have a contactomorphism $V_2\rightarrow V_2'$.

Let $f_{\partial N_1}$ be the corresponding function on homology classes of
closed curves on~$\partial N_1$. By slight abuse of notation,
we regard $\mu$ and $\lambda$ as classes both on $T$ and on~$\partial N_1$,
and we identify $M$ with $\mu\times\lambda\times [0,1]$. In the standard
neighbourhood $N_1$ (as in Case~2), we may in fact take the vector field
$v_{\partial N_1}$ in the definition of $f_{\partial N_1}$ to be
equal to~$\partial_{\theta}$; at any rate, $v_{\partial N_1}$ extends
as a nonzero section of $\xi_0$ over~$N_1$, coinciding with the
tangent vector field to $L_1$ along that Legendrian curve. It follows
that $f_{\partial N_1}(\mu )=\rot (L_1)$ and $f_{\partial N_1}(\lambda)=0$.

The Legendrian rulings
on $T$ and $\partial N_1$ define a section $\sigma$ of $\xi_0$
over~$\partial M$, with $\sigma |_T=v_T$ and $\sigma |_{\partial N_1}=
v_{\partial N_1}$. Let $e=e(\xi_0,\sigma)\in H^2(M,\partial M)$
be the corresponding relative Euler class.
Then, up to sign convention,
\[ e(\mu\times [0,1])=f_{\partial N_1}(\mu)-f_T(\mu)=\rot (L_1)-f_T(\mu)\]
and
\[ e(\lambda\times [0,1])=f_{\partial N_1}(\lambda )-f_T(\lambda)
=-f_T(\lambda ).\]

It follows that $\rot (L_1)$ and $\rot (L_2)$ determine the
relative Euler class of~$M$. The boundary slopes of $M$ are
$\tb (L_1)$ and~$q/p$. Thus, as in Case~2 we conclude that
$M$ and $M'$ are contactomorphic. 
As in the proof of Lemma~\ref{lem:isotopic}, we find a contact
isotopy that moves $L_1$ to itself and sends $T$ to~$T'$.
In that lemma, $L_2'$ and the image of $L_2$ under the isotopy
were Legendrian ruling curves and hence isotopic. Here, they
are Legendrian divides, which are isotopic by~\cite[Lemma~3.17]{etho01}.
That lemma states that these divides can be realised as the intersection
of $T'$ with a pre-Lagrangian, i.e.\ linearly foliated torus 
with closed leaves (inside
a small neighbourhood of~$T'$). Thus, the two Legendrian divides of $T'$
are Legendrian isotopic via a family of leaves on the pre-Lagrangian torus.
\end{proof}

\begin{rem}
The actual choice of Legendrian rulings on the two boundary tori
is irrelevant for the computation of the relative Euler class,
since all Legendrian rulings are homotopic to one another.
\end{rem}

The following lemma is the direct analogue of Lemma~4.12 in~\cite{etho01},
and the proof goes through essentially without any changes, except
that once again we have to use the classification of tight
contact structures on a thickened torus. By $S^b_{\pm}$ we denote
the $b$--fold positive resp.\ negative stabilisation of a Legendrian
knot.

\begin{lem}
\label{lem:3b2i-stabilise}
Let $\LL$ and $\LL'$ be $(p,q)$--cable links with $\tb (L_2)=\tb (L_2')=pq$.
Write $|p|=aq+b$ with $a,b\in\N$ and $0<b<q$.

If $\rot (L_2')=\rot (L_2)-2b$,
then $L_1\sqcup S^b_-(L_2)$ and $L_1\sqcup S^b_+(L_2')$ are Legendrian
isotopic. 

If $\rot (L_2')=\rot (L_2)-2(q-b)$,
then $L_1\sqcup S^{q-b}_-(L_2)$ and $L_1\sqcup S^{q-b}_+(L_2')$ are
Legendrian isotopic.
\hfill\qed
\end{lem}

From this lemma, the proof of Theorem~\ref{thm:1} in Case~3~(b2-i) concludes
as in Theorem~4.13 of~\cite{etho01}. The structure of the argument will
become clearer in our discussion of the remaining cases, where we
have to provide details, since~\cite{etho01} can no longer be invoked
directly.

\vspace{2mm}

{\bf Case 3 (b2-ii):} $q\geq 2$, $p<-1$ (sic!), $mp+q\leq 0$, and $|p|<q$.
Write $q=a|p|+b$ with $a,b\in\N$ and $0<b<|p|$. Thus the slope $q/p$
of $T$ lies in the interval $(-a-1,-a)$.
(Here the set-up with $L_2\subset T$ and $\tb (L_2)=pq$ maximal
is as before.) By arguments as in the
proof of Case~3~(a), we find a solid torus $S_a\subset \Int (V_2)$
with convex boundary $T_a:=\partial S_a$ having two dividing curves
of slope~$-a$. Similarly, there is a solid torus $S_{a+1}\subset S^3\setminus
L_1$ containing $V_2$ in the interior, with convex boundary $T_{a+1}:=
\partial S_{a+1}$ having two dividing curves of slope~$-(a+1)$.

Let $f_a$ and $f_{a+1}$ be the functions on homology classes,
as discussed in the preceding case, corresponding to the
tori $T_a$ and $T_{a+1}$, respectively. In order to apply
the results of~\cite{etho01}, we regard $T$ as the boundary of~$V_1$,
with positive basis for $H_1(\partial V_1)$ given by
meridian $-\lambda$ and longitude~$-\mu$. In terms of this basis,
we have $L_2=(-q)\cdot (-\lambda )+|p|\cdot (-\mu )$.
So the considerations of \cite[pp.~88/89]{etho01} for negative torus knots
apply if we replace $p,q,\mu ,\lambda$ in their formulae by
$-q, |p|,-\lambda ,-\mu$, respectively. This yields the following:

\begin{itemize}
\item[(i)] $f_T(\lambda )=|p|-1$ or $1-|p|$.
\item[(ii)] $f_a(\mu )\in\{ 1-a,3-a,\ldots ,a-3,a-1\}$.
\item[(iii)] If $f_T(\lambda )=|p|-1$, then $f_T(\mu )=f_a(\mu )+q-a$.
\item[(iv)] If $f_T(\lambda )=1-|p|$, then $f_T(\mu )=f_a(\mu )+a-q$.
\end{itemize}

Here, however, the straightforward analogy with~\cite{etho01} ends, for
the presence of $L_1$ now imposes additional restrictions on the
allowable rotation numbers for~$L_2$. First we need the following lemma.

\begin{lem}
\label{lem:rot1}
If $f_T(\lambda )=|p|-1$, then $f_{a+1}(\mu )-f_a(\mu )=1$.
If $f_T(\lambda )=1-|p|$, then $f_{a+1}(\mu )-f_a(\mu )=-1$.
\end{lem}

\begin{proof}
The fact that $f_{a+1}(\mu )-f_a(\mu )=\pm 1$ is implicit
in~\cite[p.~92]{etho01}. Here is the argument. Make the Legendrian
ruling curves on $T_a$ and $T_{a+1}$ meridional, and let $A$ be
a meridional annulus between the two tori, with Legendrian
boundary~$\partial A$. Denote the boundary curves of $A$ by
$\gamma_a=T_a\cap A$ and $\gamma_{a+1}=T_{a+1}\cap A$.
With $\gamma$ denoting either boundary curve and $\gamma'$ its
push-off along~$A$, we have $\lk (\gamma ,\gamma ')=0$, and hence
$t_A(\gamma )=\tb (\gamma )\leq 0$. So Kanda's theorem
allows us to assume that $A$ is convex.

The dividing curves on $A$ will intersect $\gamma_a$ and $\gamma_{a+1}$
in $2a$ and $2a+2$ points, respectively. As shown in~\cite[p.~92]{etho01},
there must be $2a$ dividing curves running from $\gamma_a$ to
$\gamma_{a+1}$, and one dividing curve with both endpoints
on~$\gamma_{a+1}$, cutting off a disc-like region
(a so-called boundary-parallel dividing curve). Otherwise
there would be a boundary-parallel dividing curve for $\gamma_a$,
and Etnyre and Honda show that this would give rise
to an overtwisted disc.

By~\cite[Prop.~4.5]{hond00}, the relative Euler class of the
thickened torus between $T_a$ and $T_{a+1}$, when evaluated on~$A$,
equals $\chi (R_+)-\chi (R_-)$, i.e.\ the difference in Euler
characteristics of the positive and negative regions $R_{\pm}$ into
which the dividing curves partition~$A$. The $2a$ dividing
curves running from $\gamma_a$ to $\gamma_{a+1}$ give $a$ positive and
$a$ negative disc-like regions; the boundary parallel arc adds one
positive or negative disc-like region. It follows that
$f_{a+1}(\mu )-f_a(\mu )=\pm 1$, as claimed.

We now want to establish the sign in this equation. Consider the case
that $f_T(\lambda )=|p|-1$; the case $f_T(\lambda )=1-|p|$ is completely
analogous. Make the Legendrian ruling of $T_{a+1}$ by curves
of slope $q/p$, and let $A'$ be an annulus between such a Legendrian
ruling curve and a Legendrian divide (likewise of slope~$q/p$) on~$T$.
The linking of either boundary curve of $A'$ with its push-off along
$A'$ equals $pq$, the maximal Thurston-Bennequin invariant of
a $(p,q)$--torus knot. Thus, Kanda's formula and theorem
once again imply that $A'$ may be assumed to be convex. Let $e$ be the
relative Euler class of the thickened torus between $T$ and~$T_{a+1}$.
Then
\[ e(A')=p\bigl( f_{a+1}(\mu )-f_T(\mu )\bigr) +q\bigl( f_{a+1}(\lambda )
-f_T(\lambda )\bigr) .\]
The dividing curves of $T_{a+1}$ intersect $\lambda$ on $T_{a+1}$ in
two points. In a convex longitudinal disc with
Legendrian boundary, therefore, there will
be only one dividing curve, separating the disc into one positive and
one negative region. Since $f_{a+1}(\lambda)$ may be computed by
counting the difference between the number of positive and
negative regions on a longitudinal disc (cf.~\cite[Prop.~4.23]{hond00}),
we have $f_{a+1}(\lambda )=0$.

Arguing by contradiction, we now assume that $f_{a+1}(\mu )-f_a(\mu )=-1$.
Then, together with equation~(iii) before this lemma, we obtain
\begin{eqnarray*}
e(A') & = & p\bigl( f_{a+1}(\mu )-f_T(\mu )\bigr)-qf_T(\lambda )\\
      & = & p\bigl( f_a(\mu )-1-f_T(\mu )\bigr)-qf_T(\lambda )\\
      & = & p(a-q-1)-q(|p|-1)\\
      & = & ap+q-p\; = \; |p|+b.
\end{eqnarray*}
The dividing curves on $T$ do not intersect the boundary of~$A'$; the dividing
curves on $T_{a+1}$ intersect the boundary of $A'$ in
\[ 2\left|\begin{array}{rc}p&1\\q&-(a+1)\end{array}\right|
=2\bigl( (a+1)|p|-q\bigr) =2(|p|-b)\]
points. Thus, all (non-closed)
dividing curves on $A'$ have both endpoints on $T_{a+1}$, and
$e (A')$ can be at most equal to $|p|-b$ (if there are
$|p|-b$ positive disc-like regions and one negative annulus-like
region on~$A'$) --- unless there are closed dividing curves on~$A'$.

Since our previous computation gave $e(A')=|p|+b>|p|-b$, such closed
dividing curves would indeed have to exist. But this is impossible in
a tight contact manifold by a criterion of Giroux,
cf.~\cite[Thm.~3.5]{hond00}. This contradiction proves the lemma.
\end{proof}

\begin{cor}
If $f_T(\lambda )=|p|-1$, then
\[ \rot (L_2)=pf_{a+1}(\mu )-p-b.\]
If $f_T(\lambda )=1-|p|$, then
\[ \rot (L_2)=pf_{a+1}(\mu )+p+b.\]
\end{cor}

\begin{proof}
This follows from $\rot (L_2)=pf_T(\mu )+qf_T(\lambda)$ by
a straightforward computation, using the preceding lemma and equations
(iii), (iv).
\end{proof}

Combined with the next lemma, this shows how the range
of the rotation number $\rot (L_2)$ is restricted by the
value of $\rot (L_1)$ (or vice versa). Observe that the equations
$mp+q\leq 0$ and $q=a|p|+b$ imply $m\geq a+1$.

\begin{lem}
\label{lem:bound-rot1}
We have $|\rot (L_1)-f_{a+1}(\mu )|\leq m-a-1$.
\end{lem}

\begin{proof}
Consider a convex meridional annulus, with Legendrian boundary curves,
between the tori $\partial N_1$ and~$T_{a+1}$. Analogous to the proof
of Lemma~\ref{lem:rot1}, there must be $2(a+1)$ dividing curves
on this annulus running from one boundary component to the other,
and $(m-a-1)$ dividing curves with both endpoints on~$\partial N_1$.
By the reasoning
employed in that proof, the difference
\[ |\rot (L_1)-f_{a+1}(\mu )|=|f_{\partial N_1}(\mu )-f_{a+1}(\mu )|,\]
can be at most equal to $m-a-1$, the maximum being attained
if all the extra $(m-a-1)$ dividing curves are boundary-parallel
and separate off disc-like
regions of the same sign.
\end{proof}

Lemma~\ref{lem:rot1} and condition (ii) before that lemma
(with $a$ replaced by~$a+1$)
yield:
\begin{itemize}
\item If $f_T(\lambda )=|p|-1$, then
\[ f_{a+1}(\mu )\in\{ -a+2,-a+4,\ldots ,a-2,a\}.\]
\item If $f_T(\lambda )=1-|p|$, then
\[ f_{a+1}(\mu )\in\{ -a,-a+2,\ldots ,a-4,a-2\}.\]
\end{itemize}

This apparent dichotomy is only superficial.
By some simple arithmetic, the condition on $\rot (L_2)$ can be rewritten as
\begin{equation}
\label{eqn:rot2}
\rot (L_2)\in\left\{ \pm \bigl( -q+(2l+1)|p|\bigr) \co l\in\Z ,0\leq l<
-\frac{-q+|p|}{|p|}\right\} .
\end{equation}
Not surprisingly, this
is the same condition as in Case 3~(b2-i), with $(p,q)$
replaced by $(-q,|p|)$. However, it is no longer obvious how
to add an unknot $L_1$ to Figure~8 (bottom) of \cite{etho01} realising
any combination of $\tb (L_1)$ and $\rot (L_1)$, because now $L_1$
corresponds to a longitudinal curve in that figure. Indeed, the
preceding lemmas give the restriction
\begin{equation}
\label{eqn:rot1}
\bigl| \rot (L_1)\pm (2l-a)\bigr|\leq m-a-1,
\end{equation}
with the choice of sign as in condition~(\ref{eqn:rot2}). Notice that
by the parity condition
\[ \tb (L)+\rot (L)\equiv 1\;\;\mbox{\rm mod}\; 2,\]
satisfied by any Legendrian knot $L$ in $(S^3,\xi_0)$,
cf.~\cite[Prop.~2.3.1]{elia93}, only every other
value in the range allowed by condition~(\ref{eqn:rot1}) can actually be
attained by~$\rot (L_1)$.

We claim that
all the rotation numbers allowed by conditions (\ref{eqn:rot2})
and~(\ref{eqn:rot1}) are indeed realised by a Legendrian $(p,q)$--cable link.
By the observations above, Figure~8 (bottom) of \cite{etho01} gives
a Legendrian realisation $L_2$ of the $(p,q)$--torus knot once we
replace $(p,q)$ in that figure by $(-q,|p|)$.

Given this $L_2$, we choose a solid torus $S_{a+1}$ as at the beginning
of the discussion of the present Case 3~(b2-ii). With respect to
the solid torus $S^3\setminus\Int (S_{a+1})$, the slope of
$T_{a+1}$ is equal to $-1/(a+1)$. By the classification of tight
contact structures on solid tori~\cite[Prop.~4.23]{hond00},
the contact structure on $S^3\setminus\Int (S_{a+1})$ is unique
up to isotopy, since a meridional disc of that torus is
intersected only once by each of the two dividing curves, cf.\ the proof
of Lemma~\ref{lem:rot1}. We may therefore identify
$S^3\setminus\Int (S_{a+1})$ with the solid torus
$\{ (x,y,\theta )\co x^2+y^2\leq 1\}$ in the standard model
\[ \cos ((a+1)\theta )\, dx-\sin ((a+1)\theta )\, dy=0.\]
Let $L_1^*$ be the Legendrian unknot in $S^3\setminus\Int (S_{a+1})$
corresponding to the spine $\{ x=y=0\}$ in that model, oriented so as
to be isotopic to $\mu$ in our customary notation.
Then $\tb (L_1^*)=-(a+1)$ and $\rot (L_1^*)=f_{a+1}(\mu )$.
By $(m-a-1)$--fold stabilisation of $L_1^*$, we obtain a Legendrian
unknot $L_1$ topologically isotopic to~$L_1^*$, with
$\tb (L_1)=-m$ and $\rot (L_1)$ taking on any given value in
the allowable range (subject to the parity condition mentioned above).

Lemma~\ref{lem:3b2i-isotopy} is valid unchanged in the present case.
Lemma~\ref{lem:3b2i-stabilise} has to be replaced by the following, with
the obvious changes in the proof.

\begin{lem}
\label{lem:3b2ii-stabilise}
Let $\LL$ and $\LL'$ be $(p,q)$--cable links with $\tb (L_2)=\tb (L_2')=pq$.
Write $q=a|p|+b$ with $a,b\in\N$ and $0<b<|p|$.

If $\rot (L_2')=\rot (L_2)-2b$,
then $L_1\sqcup S^b_-(L_2)$ and $L_1\sqcup S^b_+(L_2')$ are Legendrian
isotopic. 

If $\rot (L_2')=\rot (L_2)-2(|p|-b)$,
then $L_1\sqcup S^{|p|-b}_-(L_2)$ and $L_1\sqcup S^{|p|-b}_+(L_2')$ are
Legendrian isotopic.
\hfill\qed
\end{lem}

Now, again, the proof of Theorem~\ref{thm:1} in Case 3~(b2-ii) concludes
as in Theorem~4.13 of~\cite{etho01}. The argument goes as follows.
Graph the possible pairs $(\rot (L_2),\tb (L_2))$ as dots in the
plane, with $\rot$ labelled along a horizontal axis and $\tb$
along a vertical axis. Then these dots will form a `mountain range',
cf.\ Figure~9 of~\cite{etho01}.
The peaks correspond to pairs $(\rot (L_2),\tb (L_2))$ with $\tb (L_2)$
taking on the maximal value~$pq$. By Lemma~\ref{lem:3b2i-isotopy},
there are unique Legendrian realisations $\LL$ of our cable link
with $L_2$ having the classical invariants corresponding to such a peak
(and a given allowable pair of classical invariants for the unknot~$L_1$).

By condition~(\ref{eqn:rot2}), neighbouring peaks in that mountain
range have distance $2b$ or $2(|p|-b)$. Lemma~\ref{lem:3b2ii-stabilise}
shows that the valley between any two such peaks also corresponds to a unique
Legendrian realisation~$\LL$.

Since positive and negative stabilisations commute with each other,
it is now elementary to see that any two Legendrian $(p,q)$--cable links
$\LL$, $\LL'$ with the same classical invariants can be destabilised
to links $\widetilde{\LL}$, $\widetilde{\LL}'$
with invariants that are the same for both
destabilisations and which allow a unique Legendrian realisation.
This means that $\widetilde{\LL}$ and $\widetilde{\LL}'$ are Legendrian
isotopic, and so will be $\LL$ and~$\LL'$.

\vspace{2mm}

{\bf Case 3 (b2-iii):} $q\geq 2$, $p=-1$, and $-m+q\leq 0$.
The set-up with $L_2\subset T$ realising the maximal Thurston-Bennequin
invariant $\tb (L_2)=-q$ is as described at the beginning of Case 3~(b2).
Since $T$ now has integral slope $-q$, the considerations of the previous
case give $f_T(\lambda )=0$ and
\[ f_T(\mu )\in \{ -q+1,-q+3,\ldots ,q-3,q-1\} .\]
As seen in Case 3~(b2-i), for a standard neighbourhood $N_1$ of $L_1$ we
have $f_{\partial N_1}(\mu )=\rot (L_1)$ and $f_{\partial N_1}(\lambda )=0$.
As in Lemma~\ref{lem:bound-rot1} we have
\[ |\rot (L_1)-f_T(\mu )|\leq m-q.\]
Thus, given $\rot (L_1)$, the conditions above give the allowable
range of~$f_T(\mu)$, which in turn determines $\rot (L_2)$ via
\[ \rot (L_2)=-f_T(\mu )+qf_T(\lambda )=-f_T(\mu ).\]

All the rotation numbers allowed by these conditions can be realised
by a Legendrian link $\LL$ as follows. Notice that the condition $p=-1$
implies that $L_2$ is an unknot; this has a Legendrian
realisation with $\tb (L_2)=-q$ and $\rot (L_2)$ in the described range.

We now want to find a Legendrian unknot $L_1$ with $\tb (L_1)=-m$ and
$\rot (L_1)=r$, where $r$ has to satisfy the Thurston-Bennequin
inequality $|r|\leq m-1$ and the parity condition $m+r\equiv 1$ mod~$2$.
Moreover, $L_2$ has to be in the class $-\mu +q\lambda$.

Set $k=r+\rot (L_2)$. The conditions above translate into $|k|\leq m-q$,
and together with the parity condition $q+\rot (L_2)\equiv 1$ mod~$2$
we get $k\equiv m-q$ mod~$2$. We can therefore find unique non-negative
integers $n_1,n_2$ satisfying
\[ n_1+n_2=m-q\;\;\mbox{\rm and}\;\; n_1-n_2=k.\]
Let $L_1^0$ be the Legendrian push-off of $L_2$ (defined by a parallel
copy of $L_2$ in the front projection picture as in the example
at the end of case 3~(b1)), but with {\em reversed} orientation.
Then $\lk (L_1^0,L_2)=-\tb (L_2)=q$, so $L_2$ may be regarded as a cable
of $L_1^0$ representing the class $-\mu +q\lambda$. Now define
$L_1$ as the stabilisation $L_1=S^{n_1}_+S^{n_2}_-(L_1^0)$; this
does not change the topological picture.

We compute
\[ \tb (L_1)=\tb (L_1^0)-n_1-n_2=\tb (L_2)-m+q=-m\]
and
\[ \rot (L_1)=\rot (L_1^0)+n_1-n_2=-\rot (L_2)+k=r.\]

Now we come to the proof of Theorem~\ref{thm:1} in the present case.
Since $\rot (L_2)$ determines $f_T(\mu )$, the proof of
Lemma~\ref{lem:3b2i-isotopy} goes through as before. The analogue
of Lemmas \ref{lem:3b2i-stabilise} and \ref{lem:3b2ii-stabilise} is
the following.

\begin{lem}
Let $\LL$ and $\LL'$ be $(-1,q)$--cable links with $\tb (L_2)=\tb (L_2')
=-q$. If $\rot (L_2')=\rot (L_2)-2$, then $L_1\sqcup S_-(L_2)$ and
$L_1\sqcup S_+(L_2')$ are Legendrian isotopic.
\end{lem}

\begin{proof}
Observe that the assumption $\rot (L_2')=\rot (L_2)-2$ implies that
$\rot (L_2)$ is not equal to the minimal possible value $1-q$,
hence $f_T(\mu )<q-1$.

In Section~4.4.5 of \cite{hond00} it is explained how the classification
of contact structures on a solid torus is achieved by decomposing
a thickened torus $T^2\times [0,1]$ into `positive' and `negative layers'
of negative integer boundary slopes.
The sign of these layers is exactly determined by the function $f_T(\mu )$.
Moreover, positive and negative layers can be shuffled around.
The fact that $f_T(\mu )<q-1$ means that there must be at least one negative
layer in $V_2$, which we may take to be the outermost one. This
implies that we can find a convex torus $T_0$ in $V_2$ with two
dividing curves of slope $-(q-1)$ and with $f_{T_0}(\mu )=f_T(\mu )+1$.

Now consider an annulus $A$ between $T$ and $T_0$ of slope $-q$,
with $A\cap T=L_2$ and $A\cap T_0$ equal to a Legendrian ruling curve.
Both boundary curves satisfy the condition $t_A\leq 0$ that allows us
to assume that $A$ is convex: the Legendrian ruling curve must satisfy
this condition with respect to the convex torus $T_0$, and hence also with
respect to~$A$; the Legendrian divide $L_2$ actually has $t_A(L_2)=0$
by Kanda's formula.

The dividing curves of $A$ do not intersect the boundary component~$L_2$,
and they intersect the boundary component on $T_0$ in
\[ 2\left|\begin{array}{rc}1&1\\-q&-(q-1)\end{array}\right| =2\]
points. So there is one boundary-parallel dividing curve along
$A\cap T_0$. This implies that $A\cap T_0$ is a negative
stabilisation of~$L_2$: the sign of the stabilisation follows from
the condition $f_{T_0}(\mu )=f_T(\mu )+1$ and the
fact that the Legendrian knots in question
lie in the class $-\mu +q\lambda$.

Similarly, we find a convex torus $T_0'$ in $V_2'$ with two
dividing curves of slope $-(q-1)$ and
\[ f_{T_0'}(\mu )=f_{T'}(\mu )-1=-r+1=f_T(\mu )+1=f_{T_0}(\mu ).\]
We then find an annulus $A'$ with $A'\cap T_0'=S_+(L_2')$.

As in the proof of Lemma~\ref{lem:3b2i-isotopy}, 
the condition $f_{T_0'}(\mu )=f_{T_0}(\mu )$ allows us to find a contact
isotopy moving $L_1$ to itself and $T_0$ to $T_0'$.
Now argue as at the end of the proof of Lemma~\ref{lem:isotopic} to
obtain the contact isotopy sending $L_1\sqcup S_-(L_2)$ to $L_1\sqcup
S_+(L_2')$.
\end{proof}

The proof of Theorem~\ref{thm:1} now concludes exactly as in the
previous two subcases.

\section{Contact geometry of $J^1(S^1)$}
\label{section:jet}
In this section we discuss the classical invariants of Legendrian knots
in the $1$--jet space $(J^1(S^1),\xi_1)$ and prove
Theorem~\ref{thm:2}. The definition of these
invariants is based on Proposition~\ref{prop:jet}.

\begin{proof}[Proof of Proposition~\ref{prop:jet}]
Define a map $f\co J^1(S^1)\rightarrow S^3$ by
\[ f(q,p,z)=\lambda\bigl(\frac{p}{2}\sin q-z\cos q ,\cos q, \sin q,
\frac{p}{2}\cos q+z\sin q\bigr)\]
with $\lambda =1/\sqrt{1+p^2/4+z^2}$.
This $f$ is a diffeomorphism of $J^1(S^1)$ onto $S^3\setminus K_0$, where
\[ K_0=\{ (x_1,y_1,x_2,y_2)\in S^3\co y_1=x_2=0\} .\]
A straightforward calculation yields
\[ f^*(x_1\, dy_1 -y_1\, dx_1+x_2\, dy_2-y_2\, dx_2)=\lambda^2
(dz-p\, dq),\]
so $f$ is indeed a contactomorphism from $(J^1(S^1),\xi_1)$
to $(S^3\setminus K_0,\xi_0)$ preserving (co-)orientations.

\begin{rem}
The idea behind the definition of $f$ is the following: 
If we change coordinates on $J^1(S^1)$ from $(q,p,z)$ to $(q,x,y)$ with
\[ p=-2x\sin (2q)-2y\cos (2q),\;\;\; z=x\cos (2q)-y\sin (2q),\]
then $dz-p\, dq=\cos (2q)\, dx-\sin (2q)\, dy$,
which exhibits an obvious Legendrian $S^1$--fibration of $J^1(S^1)$
over $\R^2$. Likewise,
there is a Legendrian Hopf fibration of $(S^3,\xi_0)$ over~$S^2$. Take
$K_0\subset S^3$ to be a fibre, and define $f$ sending fibres to
fibres. In terms of the coordinates $(q,x,y)$, the map $f$ is given by
\[ (q,x,y)\longmapsto \frac{1}{\sqrt{1+x^2+y^2}}\,
(-x\cos q+y\sin q, \cos q,\sin q,-x\sin q-y\cos q).\]
\end{rem}

It remains to check that $\tb (K_0)=-1$. Recall that $\tb (K_0)$ is the
linking number of $K_0$ with a parallel copy $K_0'$ obtained by pushing
$K_0$ in the direction of a vector field along $K_0$ transverse to~$\xi_0$.
A candidate for such a vector field is $x_1\partial_{y_1}-y_2\partial_{x_2}$,
so we may take
\[ K_0'=\{ (x_1,x_1,-y_2,y_2)\in S^3\} .\]
The linking number $\tb (K_0)=\lk (K_0,K_0')$ can be computed as
the intersection number of Seifert surfaces $\Sigma_0$ and
$\Sigma_0'$ of $K_0$ and $K_0'$ resp.\ in
the $4$--ball~$D^4$. For $\Sigma_0$ we may take the disc
\[ \Sigma_0=\{ (x_1,0,0,y_2)\in D^4\},\]
for $\Sigma_0'$, the disc
\[ \Sigma_0'=\{ (x_1,x_1,-y_2,y_2)\in D^4\}.\]
These discs intersect at the origin $0\in D^4$
only, and one checks that a positively
oriented basis for the tangent space $T_0\Sigma_0$ followed by a positive
basis for $T_0\Sigma_0'$ is a negative basis for $T_0\R^4$.
\end{proof}

By identifying $S^1$ with $\R /\Z$, we can visualise
a Legendrian knot $K\subset J^1(S^1)$ in its front projection to
a strip $[0,1]\times\R$ in the $qz$--plane. The usual definition for
the Thurston-Bennequin invariant of $K$ is
\[ \tb (K)=\mbox{\rm writhe}(K)-\frac{1}{2}\# (\mbox{\rm cusps}(K)),\]
where the quantities on the right are computed from the front projection
of~$K$. This is the signed number of crossing changes required to
`unlink' $K$ from its push-off $K'$ in the $z$--direction (i.e.\
transverse to~$\xi_1$), that is, the number of crossing changes that
will allow one to separate the two knots in $J^1(S^1)$.

The images of the curves $q\mapsto (q,p_0,z_0)\in J^1(S^1)$,
$q\in S^1$, under the contactomorphism $f$ are Hopf fibres, any two of
which have linking number~$-1$. This accounts for the correction term
$n^2$ between $\tb (K)$ and $\tb (f(K))$.

Regarding the rotation number, we observe the following.
The contact structure $\xi_1$ admits the global sections
$\partial_q+p\partial_z$ and~$\partial_p$.
For an oriented Legendrian knot $K\subset (J^1(S^1),\xi_1)$, the rotation
number $\rot (K)$ counts the rotations of a positive tangent vector to $K$
relative to this ordered basis of sections
as we go once around~$K$. Similarly, in $(S^3,\xi_0)$
one counts rotations relative to a trivialisation and
orientation of $\xi_0$ given by
the vector fields
\[ e_1=(-y_2,-x_2,y_1,x_1)\;\;\mbox{\rm and}\;\;e_2= (x_2,-y_2,-x_1,y_1),\]
where these vector fields are written in terms of the basis
$\{ \partial_{x_1},\partial_{y_1},\partial_{x_2},\partial_{y_2}\}$.
We compute
\begin{eqnarray*}
f_*(\partial_p) & = & -\frac{p}{4}\lambda^3\bigl(\frac{p}{2}\sin q-z\cos q,
                      \cos q,\sin q,\frac{p}{2}\cos q+z\sin q\bigr)+\\
                &   & \mbox{}+ \lambda\bigl(\frac{1}{2}\sin q,0,0,\frac{1}{2}
                        \cos q\bigr)\\
                & = & \lambda^3\bigl(\frac{p}{4}z\cos q
                      +\frac{1}{2}(1+z^2)\sin q,-\frac{p}{4}\cos q,\\
                &   & \hphantom{xxxxx} -\frac{p}{4}\sin q,
                      -\frac{p}{4}z\sin q+\frac{1}{2}(1+z^2)\cos q\bigr)\\
                & = & \frac{\lambda^2}{2}(e_2-ze_1)
\end{eqnarray*}
which is homotopic through nonvanishing sections of $\xi_0$ to the vector
field~$e_2$. This proves that $\rot (K)=\rot (f(K))$.
The same considerations as for the standard contact structure on
$\R^3$ (defined likewise by $dz-p\, dq=0$, now with $q\in\R$) apply to show
that in the front projection to the $qz$--plane we have
\[ \rot (K)=\frac{1}{2}(c_- - c_+),\]
with $c_{\pm}$ the number of cusps oriented upwards or downwards,
respectively.

\begin{proof}[Proof of Theorem~\ref{thm:2}]
Let $K$ and $K'$ be two oriented Legendrian torus knots in $J^1(S^1)$ which
have the same oriented knot type and classical invariants. With
$f\co J^1(S^1)\rightarrow S^3$ denoting the contactomorphism defined in the
proof of Proposition~\ref{prop:jet}, we have the Legendrian cable
links $\LL :=K_0\sqcup f(K)$ and $\LL':=K_0\sqcup f(K')$ in
$(S^3,\xi_0)$. The formulae relating the classical invariants
in $J^1(S^1)$ and $S^3$ imply that $\LL$ and $\LL'$ have the same invariants
and are therefore isotopic by Theorem~\ref{thm:1}. Notice that
because of $\tb (K_0)=-1$, the case discussed after the proof of
Lemma~\ref{lem:tb-estimate}, where we might have had to reverse the
roles of $K_0$ and $f(K)$, does not occur.

The proof of Theorem~\ref{thm:1} does not, however, guarantee the
existence of an isotopy fixing $K_0$. Instead, we may argue as follows.
What the proof does give is a contactomorphism of $(S^3,\xi_0)$
sending $\LL$ to $\LL'$ and equal to the identity map in a neighbourhood
of~$K_0$. This translates into a contactomorphism of $(J^1(S^1),\xi_1)$
sending $K$ to $K'$, equal to the identity map outside a compact set.

If we change coordinates on $J^1(S^1)$ from $(q,p,z)$ to $(q,u,v)$
with
\[ p=-u\sin q-v\cos q,\;\;\; z=u\cos q-v\sin q,\]
then $dz-p\, dq=\cos q\, du-\sin q\, dv$. This identifies
$(J^1(S^1),\xi_1)$ as the space of cooriented contact elements
of~$\R^2$.

A contactomorphism of that space, fixed outside a compact set,
may be regarded as a contactomorphism of the space of
contact elements of the closed $2$--disc~$D^2$, fixed near
the boundary $S^1\times \partial D^2$.
According to~\cite[Thm.~1]{giro01}, the isotopy classification of
such contactomorphisms coincides with the isotopy classification
of diffeomorphisms of $D^2$, fixed near the boundary.

By a beautiful argument of Smale~\cite{smal59}, any diffeomorphism
of $D^2$ fixed near the boundary is isotopic to the identity map.
\end{proof}

\section{Transverse knots and links}

The argument of \cite[Thm.~2.10]{etho01} concerning
the classification of transverse knots in a contact manifold
(i.e.\ knots everywhere transverse to the contact structure)
via that of Legendrian
knots carries over without changes and yields the following corollary
to Theorems \ref{thm:1} and~\ref{thm:2}. Notice that transverse knots
and links carry a natural orientation making them positively transverse
to the given (cooriented) contact structure.

\begin{cor}
Transverse cable links in $(S^3,\xi_0)$
are classified up to transverse isotopy
by their oriented link type and the self-linking numbers of the
two components.

Transverse torus knots in $(J^1,\xi_1)$
are classified up to transverse
isotopy by their oriented knot
type and self-linking number.\hfill\qed
\end{cor}

\begin{ack}
F.~D.\ is partially supported by grant no.\ 10201003 of the National
Natural Science Foundation of China. H.~G.\ is partially supported by
grant no.\ GE 1245/1-1 of the Deutsche Forschungsgemeinschaft within
the framework of the Schwerpunktprogramm 1154 ``Globale
Differentialgeometrie''. Most of this research was carried out during a visit
by F.~D.\ to the Mathematical Institute of the Universit\"at zu K\"oln.
\end{ack}


\end{document}